\documentclass[a4paper,twoside,12pt]{amsart} 
 \usepackage[T1]{fontenc}            
 \usepackage{epsf}
 \usepackage{amscd}                  
 \usepackage{xypic}                  
 \usepackage{amssymb}
\begin{document}

\pagenumbering{arabic}
  \renewcommand{\datename}{\textit{Dato}:} 
\def\Im{\operatorname{Im}}          
\newtheorem{thm}{Theorem}[section]
\newtheorem{kor}[thm]{Korollar}
\newtheorem{cor}[thm]{Corollary}
\newtheorem{construction}[thm]{Construction}
\newtheorem{conj}[thm]{Conjecture}
\newtheorem{lemma}[thm]{Lemma}
\newtheorem{prop}[thm]{Proposition}
\newtheorem{oppg}[thm]{Oppgave}
\theoremstyle{definition}
\newtheorem{defn}[thm]{Definition}
\newtheorem{bem}[thm]{Remark}
\newtheorem{notation}[thm]{Notation}
\newtheorem{example}[thm]{Example}
\renewcommand{\proofname}{Proof}
\let\u=\underline
\newcommand{\vcgreek}[1]{\boldsymbol{#1}} 
\def\i{^{-1}}
\def\RR{{\mathbb R}}
\def\kar{\operatorname{char}}
\def\coker{\operatorname{coker}}
\def\mod{\operatorname{mod}}
\def\dim{\operatorname{dim}}
\def\ker{\operatorname{ker}}
\def\deg{\operatorname{deg}}
\def\konstant{\operatorname{konstant}}
\def\ann{\operatorname{ann}}
\def\lin{\operatorname{lin}}
\def\Sing{\operatorname{Sing}}
\def\min{\operatorname{min}}
\def\H{\operatorname{H}}
\def\depth{\operatorname{depth}}
\def\pd{\operatorname{pd}}
\def\im{\operatorname{Im}}
\def\rank{\operatorname{rank}}
\def\Sec{\operatorname{Sec}}
\def\Res{\operatorname{Res}}
\def\|{\mid}
\def\CC{{\mathbb C}}
\def\GG{{\mathbb G}}
\def\ZZ{{\mathbb Z}}
\def\DD{{\mathbb D}}
\def\NN{{\mathbb N}}
\def\QQ{{\mathbb Q}}
\def\VV{{\mathbb V}}
\def\PP{{\mathbb P}}
\def\CM{Cohen-Macaulay }
\def\FF{{\mathbb F}}
\def\AA{{\mathbb A}}
\def\D{{\mathcal D}}
\def\A{{\mathcal A}}
\def\F{{\mathcal F}}
\def\J{{\mathcal J}}
\def\G{{\mathcal G}}
\def\M{{\mathcal M}}
\def\T{{\mathcal T}}
\def\N{{\mathcal N}}
\def\O{{\mathcal O}}
\def\I{{\mathcal I}}
\def\Q{{\mathcal Q}}
\def\E{{\mathcal E}}
\def\K{{\mathcal K}}
\def\C{{\mathcal C}}
\def\Z{{\mathcal Z}}
\def\V{{\mathcal V}}
\def\B{{\mathcal B}}
\def\L{{\mathcal L}}
\def\iso{\cong}
\def\congr{\equiv}
\def\sub{\subseteq}
\def\subne{\subsetneqq}
\def\hpil{\longrightarrow}
\def\surj{\hpil\mspace{-26.0mu}\hpil}
\def\Aut{\operatorname{Aut}}
\def\id{\operatorname{id}}
\def\Der{\operatorname{Der}}
\def\Tor{\operatorname{Tor}}
\def\Ext{\operatorname{Ext}}
\def\red{\operatorname{red}}
\def\spec{\operatorname{Spec}}
\def\Proj{\operatorname{Proj}}
\def\Hom{\operatorname{Hom}}
\def\End{\operatorname{End}}
\def\Hilb{\operatorname{Hilb}}
\def\Grass{\operatorname{Grass}}
\def\Pic{\operatorname{Pic}}
\def\Supp{\operatorname{Supp}}
\def\Sym{\operatorname{Sym}}
\def\GL{\operatorname{GL}}
\def\SL{\operatorname{SL}}
\def\Pic{\operatorname{Pic}}
\def\codim{\operatorname{codim}}
\def\nil{\operatorname{nil}}
\def\dh{\operatorname{dh}}
\def\det{\operatorname{det}}
\hfuzz5pc 

\title {Monomial Multiple Structures}

\author { Jon Eivind Vatne }

\address{ Matematisk institutt\\
          Johs. Brunsgt. 12 \\
          N-5008 Bergen \\
          Norway}   
        
\email{ jonev@mi.uib.no}

\begin{abstract}  In this paper we study monomial multiple structures on a
linear subspace of codimension two in
projective space.  We show that these structures determine smooth
points in their respective Hilbert schemes, with (smooth)
neighbourhoods of two such points intersecting if their Hilbert
functions are equal.  We generalize a construction for multiple
structures on points in the plane to this setting, giving a kind of
product of monomial multiple structures.  For points, this
construction can be found in Nakajima's book \cite{Nak99}.  The tools
we use for studying multiple structures are developed in \cite{Vat01}
and \cite{Vat02a}.
\end{abstract}

\maketitle

\section{Introduction} 
  The monomial multiple structures are among the best understood
  multiple structures.  They can be visualized using Young diagrams,
  their invariants are easily calculated, and their ``inductive
  construction'' (see below) is quite simple.  In this paper we focus
  on monomial multiple structures on linear subspaces of codimension
  two in projective space satisfying the Cohen-Macaulay property.
  Cohen-Macaulay always means {\em locally} Cohen-Macaulay.\\

  Let us first summarize the method for constructing multiple
  structures from \cite{Vat01} and \cite{Vat02a}, adapted to our
  special situation:  Let $X=\PP^n\subset \PP^{n+m}$ be a linear
  subspace, $X^{(i)}\subset \PP^{n+m}$ the i'th infinitesimal
  neighbourhood of $X$, and $Y$ a Cohen-Macaulay multiple structure
  with $Y_{\red}=X$, whose ideal is generated by monomials.  Then
  there is a filtration of $Y$

\[X=Y_0\subset Y_1=Y\cap X^{(1)}\subset \cdots \subset Y_{k-1}=Y\cap
  X^{(k-1)} \subset Y_k=Y\cap X^{(k)}=Y\]
for some $k$.  Since $Y$ is a monomial Cohen-Macaulay structure, the
same is true for each term in this filtration, but in general the
terms will not have such good properties.  To keep the terminology the same as
in \cite{Vat01} and \cite{Vat02a} we will therefore refer to this
filtration as the $S_1$-{\em filtration} of $Y$.\\

 Let $\I_i$ be the Ideal of $Y_i$.  Then there are two
short exact sequences
\begin{equation}
0\rightarrow \I_{j+1}/\I_X\I_j\rightarrow \I_j/\I_X\I_j\rightarrow
\L_j \rightarrow 0
\end{equation}

and

\begin{equation}
\label{strukknippl}
0\rightarrow \L_j\rightarrow \O_{Y_{j+1}}\rightarrow
\O_{Y_j}\rightarrow 0.
\end{equation}

Here the $\O_X$-Module $\L_j$ is defined by the first exact sequence,
and the $\O_{Y_i}$ are $\O_X$-Modules by restricting a projection onto
$X$ to each $Y_i$.  The Hilbert polynomial of $Y$ can be calculated as

\[\Hilb(Y)=\Hilb(X)+\sum_j \Hilb(\L_j).\]

In general, the terms $\L_j$ will be torsion free, but in the monomial
case they are actually {\em locally free}.\\

The paper is organized as follows:  In Section \ref{Ydami} we give the
basic description of monomial multiple structures in terms of Young
diagrams.  In Section \ref{PiHscfYd} we show that the points in the
Hilbert scheme determined by monomial multiple structures in
codimension two are smooth, and relate neighbourhoods of these points
for different multiple structures.  In Section \ref{SoYd,poms} we
describe a procedure for taking a kind of product of two monomial
multiple structures in codimension two, corresponding to a natural
operation on Young diagrams.\\

{\bf Acknowledgements} This article is based on parts of my Doctorate
thesis \cite{Vat01}.  I would like to thank the group in algebraic
geometry at the University of Bergen, G.Fløystad, T.Johnsen,
A.L.Knutsen, S.-A.Strømme, J.-M.Økland an most of all my advisor
A.Holme.  I will also express my gratitude to the evaluation committee
for many useful comments and suggestions: T.Johnsen, K.Ranestad and
F.-O.Schreyer.\\

\section{Young diagrams and monomial ideals}
\label{Ydami}
 \begin{defn}
An $n$-dimensional {\em Young diagram} is a finite array of unit boxes
in the first hyperquadrant in Euclidean $n$-space, with corners in the
integer grid, such that wherever there is a box, there is always a box
in the place immediately below it in each direction.  When we speak of
two-dimensional Young diagrams we will suppress the number ``two''
and speak simply of Young diagrams, and we will use terms like ``up,
down, left and right''.  The {\em weight} of a box is the sum of the
coordinates of the innermost corner in the diagram.
\end{defn}

It is well known that Young diagrams correspond to partitions: The
parts of the partition are the number of boxes in the row
corresponding to the place of the part in the partition, where the
partition is ordered by (weakly) decreasing magnitude.  As an example,
to the partition $(4,4,3,1)$ of $12$ corresponds the Young diagram
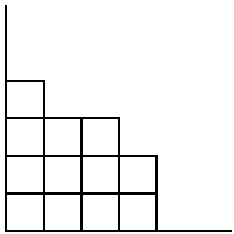
\setlength{\unitlength}{2.5cm}
\begin{picture}(4,1.5)

    \put(0,0){\line(1,0){1.2}}
    \put(0,0.2){\line(1,0){0.8}}
    \put(0,0.4){\line(1,0){0.8}}
    \put(0,0.6){\line(1,0){0.6}}
    \put(0,0.8){\line(1,0){0.2}}

    \put(0,0){\line(0,1){1.2}}
    \put(0.2,0){\line(0,1){0.8}}
    \put(0.4,0){\line(0,1){0.6}}
    \put(0.6,0){\line(0,1){0.6}}
    \put(0.8,0){\line(0,1){0.4}}

\end{picture}

  We number boxes as follows: the box $(r_1,...,r_m)$ of a Young diagram
is the box with this integer point as its lower corner (lower in all
directions).  For instance, 

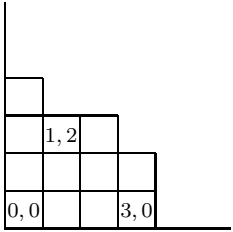
\begin{picture}(4,1.5)

    \put(0,0){\line(1,0){1.2}}
    \put(0,0.2){\line(1,0){0.8}}
    \put(0,0.4){\line(1,0){0.8}}
    \put(0,0.6){\line(1,0){0.6}}
    \put(0,0.8){\line(1,0){0.2}}

    \put(0,0){\line(0,1){1.2}}
    \put(0.2,0){\line(0,1){0.8}}
    \put(0.4,0){\line(0,1){0.6}}
    \put(0.6,0){\line(0,1){0.6}}
    \put(0.8,0){\line(0,1){0.4}}

    \put(0.1,0.04){\makebox(0,0)[b]{{\tiny $0,0$}}}
    \put(0.3,0.44){\makebox(0,0)[b]{{\tiny $1,2$}}}
    \put(0.7,0.04){\makebox(0,0)[b]{{\tiny $3,0$}}}
\end{picture}\\

Young diagrams arise from
Cohen-Macaulay monomial ideals whose radical is the ideal of a linear
subspace.  We will constantly use the following elementary lemma
without notice:

\begin{lemma}
\label{s1giracm}
Let $I_X=(x_1,\dots,x_m)\subset k[x_1,\dots\,x_m,z_0,\dots,z_n]$ be the
(saturated) ideal of a linear subspace $X=\PP^n\subset \PP^{n+m}$.
Let $J=I_Y$ be a (saturated) monomial ideal whose radical is $I$, i.e.
$Y_{\red}=X$.  Then the following are equivalent:
\begin{itemize}
\item $Y$ satisfies $(S_1)$.
\item $Y$ is (locally) Cohen-Macaulay.
\item $Y$ is arithmetically Cohen-Macaualay.
\end{itemize}
\end{lemma}
\begin{proof}
It is obviously enough to show that $Y$ is arithmetically
Cohen-Macaulay if $Y$ satisfies $(S_1)$.  So assume this.  Now each
member of a minimal monomial generator set of $J$ looks like
$x_1^{a_1}\cdots x_m^{a_m}z_0^{b_0}\cdots z_n^{b_n}$ for some
non-negative integers $a_i$ and $b_j$.  If any $b_j$ is non-zero $z_j$
is necessarily contained in an associated prime ideal of $Y$, but by
assumption $Y$ does not have any embbeded primes.  Thus all $b_j=0$
and $z_0,\dots,z_n$ is a {\em global} regular sequence on $Y$ and
$Y$ is arithmetically Cohen-Macaulay. 
\end{proof}

For the construction of associated Young diagrams, we start with
$X=\PP^n\subset \PP^{n+2}$, with ideal $I_X=(x,y)\subset 
S:=k[x,y,z_0,...,z_n]$, and we consider (Cohen-Macaulay) monomial
ideals $J\subset S$
whose radical is $I_X$.  To such an ideal $J$ we associate a Young
diagram.  First of all we will think of $x,y$ as being directions or
axes in the Euclidean plane where we draw our Young diagrams, $x$
pointing to the right and $y$ pointing up, as usual.  Then each
monomial in $x,y$ gives an integer point in the first quadrant,
$x^ry^s\leftrightarrow (r,s)$.  We write the monomial in that
position: 

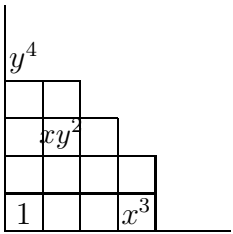
\begin{picture}(4,1.5)

    \put(0,0){\line(1,0){1.2}}
    \put(0,0.2){\line(1,0){0.8}}
    \put(0,0.4){\line(1,0){0.8}}
    \put(0,0.6){\line(1,0){0.6}}
    \put(0,0.8){\line(1,0){0.4}}

    \put(0,0){\line(0,1){1.2}}
    \put(0.2,0){\line(0,1){0.8}}
    \put(0.4,0){\line(0,1){0.8}}
    \put(0.6,0){\line(0,1){0.6}}
    \put(0.8,0){\line(0,1){0.4}}

    \put(0.1,0.04){\makebox(0,0)[b]{$1$}}
    \put(0.3,0.44){\makebox(0,0)[b]{{\small $xy^2$}}}
    \put(0.7,0.04){\makebox(0,0)[b]{$x^3$}}
    \put(0.1,0.84){\makebox(0,0)[b]{$y^4$}}

\end{picture}\\
\noindent
Now, consider the monomial ideal $J$.  Since it is an ideal, whenever
a monomial $x^ry^s$ lies in $J$, all
multiples $x^{r+u}y^{s+v}$ also lie in $J$. Saying that the radical of
$J$ is $(x,y)$ is the same as saying that there are only finitely many
pairs $(r,s)$ such that $x^ry^s\notin J$.  Put together, the monomials
not in $J$ form a Young diagram.  For instance,
$J=(x^4,x^3y^2,x^2y^3,y^4)$ gives the Young diagram (where we write
the generators of the ideal in their respective positions outside the
Young diagram)

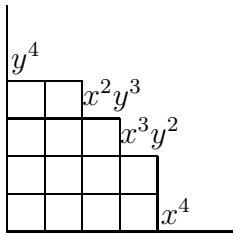
\begin{picture}(4,1.5)

    \put(0,0){\line(1,0){1.2}}
    \put(0,0.2){\line(1,0){0.8}}
    \put(0,0.4){\line(1,0){0.8}}
    \put(0,0.6){\line(1,0){0.6}}
    \put(0,0.8){\line(1,0){0.4}}

    \put(0,0){\line(0,1){1.2}}
    \put(0.2,0){\line(0,1){0.8}}
    \put(0.4,0){\line(0,1){0.8}}
    \put(0.6,0){\line(0,1){0.6}}
    \put(0.8,0){\line(0,1){0.4}}

    \put(0.9,0.04){\makebox(0,0)[b]{$x^4$}}
    \put(0.76,0.44){\makebox(0,0)[b]{$x^3y^2$}}
    \put(0.56,0.64){\makebox(0,0)[b]{$x^2y^3$}}
    \put(0.1,0.84){\makebox(0,0)[b]{$y^4$}}
\end{picture}\\
\noindent
Obviously this is a bijective correspondence between Cohen-Macaulay
monomial ideals in two variables and Young diagrams.  Similarly,

\begin{prop} There is a bijective correspondence between Cohen-Macaulay
monomial ideals in $l$ variables having a given linear subspace of
codimension l as support, and $l$-dimensional Young diagrams.
Under this correspondence, the number of boxes in the Young diagram
is the multiplicity of the scheme defined by the corresponding ideal.
The correspondence is inclusion-reversing.
\end{prop}
\begin{proof}
The correspondence is as above.  The equality between the number of
variables occuring in the ideal and the codimension ensures
that there are only finitely many boxes.\\

Each box in the Young diagram gives a monomial not in $J$.  These
monomials form a basis for the homogeneous coordinate ring $S/J$ as an
$S/I_X$-module by Theorem \ref{multiplikativ}.  The number of basis elements is equal to the rank of
this module, or equivalently, the rank of the structure sheaf as a
Module over the structure sheaf of the reduced subscheme.  This rank
is the multiplicity.\\

If one Young diagram contains another, the corresponding ideal has
fewer monomials than the second one's, and vice versa.
\end{proof}

The statement about multiplicities in the proposition will be sharpened in the
following sections; we will see, that the entire Hilbert polynomial
and Hilbert function can be read off the diagram.  \\ 

Some standard operations on ideals are easily described using the
associated Young diagrams:

\begin{prop}
Let $I$, $J$ be Cohen-Macaulay monomial ideals in $x_1,...,x_m$ with
radical $(x_1,\dots,x_m)$, with associated Young
diagrams $T$ and $S$.
\begin{itemize}
\item{}The Young diagram of $I+J$ is $T\cap S$.
\item{}The Young diagram of $I\cap J$ is $T\cup S$.
\item{}The Young diagram of the ideal of the $k$th infinitesimal
  neighbourhood is given by all boxes on or below the $k$th diagonal.
\item{}The Young diagram of the $k$th part of the $S_1$-filtration of
  $I$,  $I_k$, is given by the part of the Young diagram of $I$ on or
  below the $k$th diagonal.
\end{itemize}
\end{prop}

By ``the $k$th diagonal'' we mean the set of boxes $r_1,...,r_m$ such
that the weight $r_1+\dots+r_m$ is constant, equal to $k$.  In
dimension two a diagonal is thus an ordinary diagonal going from the
upper left to the lower right.  In higher dimensions it is a ``large
diagonal''.

\begin{proof}
A monomial not in $I+J$ is a monomial neither in $I$
nor in $J$, so it corresponds to a box in both $T$ and $S$.
Similarly, a monomial in $I\cap J$ is in both $I$ and $J$, so one not
in $I\cap J$ is outside $I$ or outside $J$.  Thus the corresponding
box is in either $T$ or $S$.  The third part follows by
noting that the monomials of a given degree have the same 
``distance'' from the origin, they have the same {\em weight},  so
they form a diagonal.  Thus the
Young diagram of an infinitesimal neighbourhood is given by all
boxes on or below a diagonal.  The last statement is immediate from
the first and the third.
\end{proof}

\begin{example} Let $I=(x^2,xy^2,y^3)$ and $J=(x^4,xy,y^2)$, so that
$I+J=(x^2,xy,y^2)$ and $I\cap J=(x^4,x^2y,xy^2,y^3)$.  First we will
consider $I+J$:

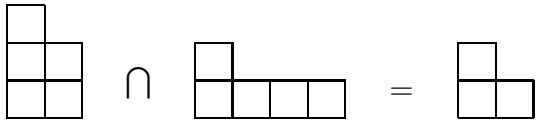
\begin{picture}(4,1)
    \put(0,0){\line(1,0){0.4}}
    \put(0,0.2){\line(1,0){0.4}}
    \put(0,0.4){\line(1,0){0.4}}
    \put(0,0.6){\line(1,0){0.2}}

    \put(0,0){\line(0,1){0.6}}
    \put(0.2,0){\line(0,1){0.6}}
    \put(0.4,0){\line(0,1){0.4}}

    \put(0.7,0.1){\makebox(0,0)[b]{$\bigcap$}}

    \put(1,0){\line(1,0){0.8}}
    \put(1,0.2){\line(1,0){0.8}}
    \put(1,0.4){\line(1,0){0.2}}

    \put(1,0){\line(0,1){0.4}}
    \put(1.2,0){\line(0,1){0.4}}
    \put(1.4,0){\line(0,1){0.2}}
    \put(1.6,0){\line(0,1){0.2}}
    \put(1.8,0){\line(0,1){0.2}}

    \put(2.1,0.1){\makebox(0,0)[b]{$=$}}

    \put(2.4,0){\line(1,0){0.4}}
    \put(2.4,0.2){\line(1,0){0.4}}
    \put(2.4,0.4){\line(1,0){0.2}}

    \put(2.4,0){\line(0,1){0.4}}
    \put(2.6,0){\line(0,1){0.4}}
    \put(2.8,0){\line(0,1){0.2}}
\end{picture}\\

Then we consider $I\cap J$:

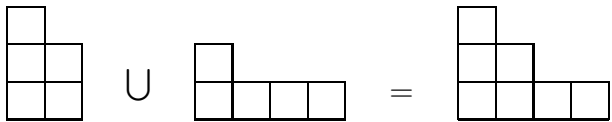
\begin{picture}(4,1)
    \put(0,0){\line(1,0){0.4}}
    \put(0,0.2){\line(1,0){0.4}}
    \put(0,0.4){\line(1,0){0.4}}
    \put(0,0.6){\line(1,0){0.2}}

    \put(0,0){\line(0,1){0.6}}
    \put(0.2,0){\line(0,1){0.6}}
    \put(0.4,0){\line(0,1){0.4}}

    \put(0.7,0.1){\makebox(0,0)[b]{$\bigcup$}}

    \put(1,0){\line(1,0){0.8}}
    \put(1,0.2){\line(1,0){0.8}}
    \put(1,0.4){\line(1,0){0.2}}

    \put(1,0){\line(0,1){0.4}}
    \put(1.2,0){\line(0,1){0.4}}
    \put(1.4,0){\line(0,1){0.2}}
    \put(1.6,0){\line(0,1){0.2}}
    \put(1.8,0){\line(0,1){0.2}}

    \put(2.1,0.1){\makebox(0,0)[b]{$=$}}

    \put(2.4,0){\line(1,0){0.8}}
    \put(2.4,0.2){\line(1,0){0.8}}
    \put(2.4,0.4){\line(1,0){0.4}}
    \put(2.4,0.6){\line(1,0){0.2}}

    \put(2.4,0){\line(0,1){0.6}}
    \put(2.6,0){\line(0,1){0.6}}
    \put(2.8,0){\line(0,1){0.4}}
    \put(3.0,0){\line(0,1){0.2}}
    \put(3.2,0){\line(0,1){0.2}}

\end{picture}\\

At last we consider the $S_1$ filtration of $J$:

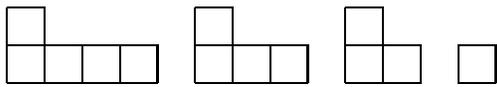
\begin{picture}(4,1)

    \put(0,0){\line(1,0){0.8}}
    \put(0,0.2){\line(1,0){0.8}}
    \put(0,0.4){\line(1,0){0.2}}

    \put(0,0){\line(0,1){0.4}}
    \put(0.2,0){\line(0,1){0.4}}
    \put(0.4,0){\line(0,1){0.2}}
    \put(0.6,0){\line(0,1){0.2}}
    \put(0.8,0){\line(0,1){0.2}}

    \put(1,0){\line(1,0){0.6}}
    \put(1,0.2){\line(1,0){0.6}}
    \put(1,0.4){\line(1,0){0.2}}

    \put(1,0){\line(0,1){0.4}}
    \put(1.2,0){\line(0,1){0.4}}
    \put(1.4,0){\line(0,1){0.2}}
    \put(1.6,0){\line(0,1){0.2}}

    \put(1.8,0){\line(1,0){0.4}}
    \put(1.8,0.2){\line(1,0){0.4}}
    \put(1.8,0.4){\line(1,0){0.2}}

    \put(1.8,0){\line(0,1){0.4}}
    \put(2,0){\line(0,1){0.4}}
    \put(2.2,0){\line(0,1){0.2}}

    \put(2.4,0){\line(1,0){0.2}}
    \put(2.4,0.2){\line(1,0){0.2}}

    \put(2.4,0){\line(0,1){0.2}}
    \put(2.6,0){\line(0,1){0.2}}

\end{picture}\\

\end{example}

\begin{prop} Given a Cohen-Macaulay monomial ideal $I$ in
    $x_1,...,x_m$ with radical $(x_1,\dots,x_m)$, consider a simple
    thickening 

\[0 \rightarrow \J/\I\I_X\rightarrow \I/\I\I_X\rightarrow
\L\rightarrow 0\]
where the map to $\L=\bigoplus \O_X(-a_i)$ is given by projection onto a
summand generated by some of the minimal monomial generators of $I$. 
The ($m$-dimensional) Young diagram of $II_X$ is constructed from the
Young diagram of 
$I$ by introducing new boxes in the inner corners.  The Young diagram of the
ideal $J$ defined by the short exact sequence is constructed by
introducing the boxes corresponding to the monomial generators of $\L$.
\end{prop}

By inner corner of the Young diagram we mean a box not in the diagram, but
such that the box immediately below it in each direction is in the
diagram.  There are also outer corners; those are the boxes such
that the box immediately below it in each direction is outside the
Young diagram, whereas the box meeting it at its corner closest to
the origin is in the diagram.

\begin{proof}
We have already used the (obvious) fact that an ideal is generated by
the monomials corresponding to inner corners of the diagram.  Multiplying a
generator with one of the generators of the ideal $I_X$ corresponds to
moving a box one step up in the direction determined by that
generator.  Thus a monomial not in $II_X$ is a monomial that corresponds
to a box at most one step up from a box in the diagram of $I$ in all
directions.  Thus the monomials not in $II_X$ correspond to monomials
in the diagram of $I$ or in the inner corners of that diagram.  The ideal $J$
contains the ideal $II_X$, but also the generators not sent to $\L$.
Thus $J$ corresponds to the diagram given by introducing boxes in the indicated
subset of the inner corners.
\end{proof}

\section{Points in Hilbert schemes coming from Young diagrams}
\label{PiHscfYd}
\begin{prop}
\label{resfromcorners}
 Given a Cohen-Macaulay monomial ideal $I$ with support a linear
  subspace of codimension two, and with (two-dimensional) Young
  diagram $T$.  Then the
  syzygies of $I$ correspond to the outer corners of $T$.  More
  precisely: if $I$ has inner corners of weight $n_{1j}$ and outer
  corners of weight $n_{2i}$, then the minimal resolution of $I$ has
  the form

\[0\rightarrow \bigoplus_i S(-n_{2i})\rightarrow \bigoplus_j
  S(-n_{1j})\rightarrow I \rightarrow 0.\]
\end{prop}

Note that there is one outer corner less than there are inner corners,
so the ranks add up correctly.
\begin{proof}
For each outer corner, there is an inner corner to its left and
another downwards.  Thus there is a relation of degree equal to the
weight of the outer corner
between the generators corresponding to these two boxes, saying that a
power of $x$ times the monomial corresponding to the inner corner to the
left is equal to a power of $y$ times the monomial corresponding to
the inner corner downwards.  It is easily verified that these
relations generate the full module of syzygies.
\end{proof}

We will need these two sets of integers later on.  We silently assume
that they are written in decreasing order, so that it makes sense to
talk about them as {\em sequences}.

\begin{cor}
\label{dimavhilb}
The dimension of the Hilbert scheme of closed subschemes in
$\PP^{n+2}$ in the point corresponding to $I$ is
\begin{eqnarray*}
\sum_{n_{2i}\geq n_{1j}}\binom{n_{2i}-n_{1j}+n+2}{n+2}
+\sum_{n_{1j}\geq n_{2i}}\binom{n_{1j}-n_{2i}+n+2}{n+2} \\
-\sum_{n_{2i}\geq n_{2j}}\binom{n_{2i}-n_{2j}+n+2}{n+2}
-\sum_{n_{1j}\geq n_{1i}}\binom{n_{1j}-n_{1i}+n+2}{n+2}+1.
\end{eqnarray*}
The Hilbert scheme is smooth in the point determined by $I$.
\end{cor}

\begin{proof}
The displayed equation is from Ellingsrud, \cite{Ell75}, where the condition is
that the ideal is of codimension two and is arithmetically
Cohen-Macaulay.  This holds in our situation because of Lemma
\ref{s1giracm}.
\end{proof}

\begin{thm}\label{multiplikativ}  Given an $m$-dimensional Young
  diagram with corresponding
  multiple scheme $Y$ with support $X$, the structure sheaf of $Y$,
  considered as an $\O_X$-Module, is
\[\O_Y\iso \bigoplus_B \O_X(-w_B).\]

Here the sum is over the boxes $B$ in the diagram, and $w_B$ is the
weight of the box $B$.  Introduce multiplicative structure on $\bigoplus_B
\O_X(-w_B)$ by introducing the monomial corresponding to $B$
as generator for $\O_X(-w_B)$.  Then the above isomorphism becomes an
isomorphism of $\O_X$-Algebras. 
\end{thm}

\begin{proof}
  The proof is by examining what happens with the Young diagram when
  we perform our standard thickenings of multiple structures.  So let
  $I$ be a monomial ideal, $I_X=(x_1,...,x_m)$ the ideal of the
  support with $\I,\I_X$ the associated sheaves.  Then our standard
  construction gives us

\[\bigoplus_{j=1}^s\O_X(-a_j)\iso \I/\I\I_X\]
where the $a_j=w_B$ are degrees of the generators of $\I$, and are thus
given by sums $r_1+\dots +r_m$ for generating monomials $x_1^{r_1}\dots
x_m^{r_m}$.  Since the
generating monomials are represented by the corners of the Young
diagram, we see that in the standard exact sequence

\[0\rightarrow \bigoplus_{j=t+1}^s\O_X(-a_j)\rightarrow \I/\I\I_X
\rightarrow \bigoplus_{j=1}^t\O_X(-a_j) \rightarrow 0\]
the ideal $J$ with $\J/\I\I_X=\bigoplus_{j=t+1}^s\O_X(-a_j)$ has Young
diagram obtained from the Young diagram of $I$ by adding the boxes
corresponding to $\bigoplus_{j=1}^t\O_X(-a_j)$.  Now this gives an
additional summand $\bigoplus\O_X(-a_j)$ to the structure sheaf by
the short exact sequence \ref{strukknippl} (which splits in this our
case).  Using this as an induction step, the theorem is proven (the
case corresponding to a single box being trivial).\\

The multiplicative structure is easily seen to be the same on both sides.

\end{proof}

\begin{bem}
The multiplicative structure on the sum $\oplus_B \O_X(-w_B)$ is a
multiplicative structure on the sum of the $\L_j$s of the
$S_1$-filtration.  A similar (but more general) construction was
considered for curves by
B\v{a}nic\v{a} and Forster, \cite{BF81} and \cite{BF86}, and also by
Manolache in higher dimensions, see \cite{Man94}.
\end{bem}

\begin{cor}
\label{computethehilbpol}
The Hilbert polynomial of the scheme associated to an $m$-dimensional
Young diagram
can be computed as follows: Define
\[b_i(d)=\binom{n+d-i}{n}=\chi (\O_{\PP^n}(d-i)).\]
  Then the Hilbert polynomial is
\[\sum_B b_{w_B}\]
the sum over all boxes $B$ in the diagram.
\end{cor}

\begin{bem}\label{comhilbfunc}  Another way of stating this corollary is to say that the
  Hilbert polynomial of the multiple structure is equal to the Hilbert
  polynomial of the sum of the invertible sheaves $\O_X(-w_B)$ from
  the theorem.  In this formulation, the statement can be strengthened
  to an equality of Hilbert {\em functions}.
\end{bem}

\begin{prop}

\label{samecomp}
All the Cohen-Macaulay monomial ideals with a given codimension two
linear subvariety as support, with a given Hilbert function, lie in the
same irreducible component of the Hilbert scheme.
\end{prop}

\begin{example} We cannot replace ``Hilbert function'' by ``Hilbert
  polynomial'' in the previous proposition.  Consider for example the
  two ideals $I=(x^5,x^4y,y^2)$ and $I'=(x^6,x^2y,xy^2,y^3)$ in
  $k[x,y,z,w]$.  The two multiple schemes defined on the line $x=y=0$
  by these ideals have the same Hilbert polynomials, but different
  Hilbert functions.  Using Corollary \ref{dimavhilb} we find that the
  dimension of the Hilbert scheme in the point corresponding to $I$ is
  38, whereas it is $39$ in the point corresponding to $I'$.
\end{example}

In order to prove this proposition, we need some combinatorial
definitions, as well as another result from \cite{Ell75}.

\begin{defn}
Define the equivalence relation $R$ on the set of Young diagrams by
$TRT'$ if and only if the ideals defined by $T$ and $T'$ have the same
Hilbert function.  By Theorem \ref{computethehilbpol} and Remark
\ref{comhilbfunc}, this is equivalent to saying that $T$ and $T'$ have
the same number of boxes in each diagonal.
\end{defn}

\begin{bem}\label{smalldiag}
For Young diagrams that are ``small'' with respect to the dimension,
this is also equivalent to having the same Hilbert polynomial.  In
fact, given two monomial ideals $I,I'$ with the same Hilbert
polynomial, but with different Hilbert functions, we can regard them
as ideals in a polynomial ring with more variables (i.e. we can
consider their projective cones).  Then, for a
sufficiently high number of variables, the Hilbert polynomials will
also be different.  Basically, this is because there are relations
between the classes of sums of line bundles in the Grothendieck group
of $\PP^N$, but these relations cannot be extended indefinitely.  Thus
Young diagrams differing by boxes corresponding to line bundles
satisfying some identities in low dimensions will not be determined by
their Hilbert polynomial in these low dimensions.
\end{bem}

\begin{defn}
Define the equivalence relation $r$ on the set of Young diagrams as follows:
By Proposition \ref{resfromcorners} each Young diagram determines the
resolution 

\[0\rightarrow \bigoplus_{i\in A} S(-n_{2i})\rightarrow
\bigoplus_{j\in B}
  S(-n_{1j})\rightarrow I \rightarrow 0\]
\noindent
for some index sets $A$ and $B$.
Consider the two sequences of integers $(\{n_{2i}\},\{n_{1j}\})$.
We say that two such pairs $(\{n_{2i}\},\{n_{1j}\})$ and
$(\{n'_{2i'}\},\{n'_{1j'}\})$ are primitively equivalent if there
exists subsets $A_0\subset A,B_0\subset B$ and a bijection
$\sigma:A_0\rightarrow B_0$
such that $n_{2i}=n_{1\sigma(i)}$ holds for all $i\in A_0$, and $\{n'_{2i'}\}=
\{n_{2i}\} \setminus A_0,\{n'_{1j'}\}=\{n_{1j}\}\setminus B_0$.  We
say that $TrT'$ if the sets of integers determined by $T$ and $T'$ are
equivalent by the equivalence relation generated by the primitive
equivalences above.  This should be compared with Ellingsrud's
``prolongement'', see \cite{Ell75}, page 424, on which it is based.
\end{defn}

What pairs of sequences of integers can appear?

\begin{lemma}
\label{charofpairs}
\begin{itemize}
\item[a)]
A pair of sequences of integers $(\{a_i\},\{b_j\})$ is associated to some
Young diagram if and only if
\begin{itemize}
\item[i)] There is exactly one more $a_i$ then there are $b_j$s.
\item[ii)] $b_i> a_i\geq a_{i+1}$ for all $i$. 
\item[iii)] $\sum_i a_i =\sum_j b_j$.
\end{itemize}
\item[b)] Any Young diagram is $r$-equivalent to a Young diagram
  without any equalities of the kind $a_i=b_j$.
\end{itemize}
\end{lemma}

In part $a) ii)$ the inequality $a_i\geq a_{i+1}$, a part of our
assumptions, is included only for emphasis.

\begin{proof}
For part $a)$, the necessity of the conditions is obvious.
Conversely, given such a
pair we can construct the Young diagram (where we write the weights on
inner and outer corners)

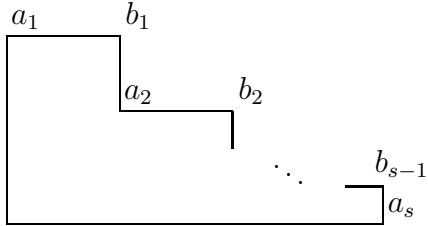
\begin{picture}(4,1.4)
    \put(0,0){\line(1,0){2}}
    \put(0,1){\line(1,0){0.6}}
    \put(0.6,0.6){\line(1,0){0.6}}
    \put(1.8,0.2){\line(1,0){0.2}}
   
    \put(0,0){\line(0,1){1}}
    \put(0.6,0.6){\line(0,1){0.4}}
    \put(1.2,0.6){\line(0,-1){0.2}}
    \put(2,0.2){\line(0,-1){0.2}}

    \put(1.4,0.2){\makebox{$\ddots$}}

    \put(0.1,1.04){\makebox(0,0)[b]{$a_1$}}
    \put(0.7,1.04){\makebox(0,0)[b]{$b_1$}}
    \put(0.7,0.64){\makebox(0,0)[b]{$a_2$}}
    \put(1.3,0.64){\makebox(0,0)[b]{$b_2$}}
    \put(2.1,0.24){\makebox(0,0)[b]{$b_{s-1}$}}
    \put(2.1,0.04){\makebox(0,0)[b]{$a_s$}}
  
\end{picture}

The box corresponding to the number $a_1$ has coordinates $(0,a_1)$.
The box corresponding to the number $b_1$ has coordinates
$(b_1-a_1,a_1)$.  The box corresponding to the number $a_2$ has
coordinates $(b_1-a_1,a_1-b_1+a_2)$.  Continuing in this way, we
finally come to the box corresponding to the number $a_s$, which has
coordinates 
\[(b_1-a_1+b_2-a_2+\dots+b_{s-1}-a_{s-1},a_1-b_1+a_2-b_2+\dots
+a_s)=(a_s,0)\]

The conditions ensure that the first coordinates of the $a_i$ (equal
to the first coordinates of the $b_{i-1}$) form an increasing sequence
from $0$ to $a_s$, whereas the second coordinates of the $a_i$ (equal
to the second coordinates of the $b_i$) form a decreasing sequence
from $a_1$ to $0$.\\

Note that the Young diagram constructed in this way is only one of
many different Young diagrams determining these sequences.\\

Part $b)$ follows from part $a)$:  just remove any pairs $a_i=b_j$ and
write the resulting sequences in order of decreasing magnitude.  The
conditions $i)-iii)$ are still verified.
\end{proof}

Now the result from \cite{Ell75}, translated to the (special) case of
monomial ideals:

\begin{prop}
\label{neigharitcm}
Given Young diagrams $T$ and $T'$.  Let $P,P'$ be the corresponding
points in the Hilbert scheme.  Then
there are smooth open neighbourhoods around $P$ and $P'$ that
intersect if and only if $TrT'$.  For a given 
Hilbert function, the union of all open sets thus determined for all
Young diagrams is a smooth, connected scheme.
\end{prop}
\begin{proof}
See \cite{Ell75} page 426.(Theorem 2.iii).  Note that we only use a
small part of this theorem.\\

In the statement of the theorem of Ellingsrud's, in order to get a
non-empty intersection, we need to know that
there is a scheme with the common refinement as the pair of sequences
arising from its minimal resolution.  This is automatic in the case of
a given Hilbert function,
since (by Lemma \ref{charofpairs}) there is a monomial ideal with the
same Hilbert function, without any equalities like $n_{2i}=n_{1j}$. 
\end{proof} 

\begin{prop}
The equivalence relations $r$ and $R$ are the same.
\end{prop}
\begin{proof}
Since the Hilbert function of $I$ can be computed from the Young
diagram by Remark \ref{comhilbfunc}, it is easily seen that $TrT'
\Rightarrow TRT'$.\\

For the other implication, note that by Remark \ref{smalldiag}, if the
Hilbert functions of two monomial ideals are different, then the zero
schemes are (succesive) hyperplane sections of monomial schemes
with different Hilbert polynomials.  These will nevertheless have the
same Young diagrams as the original monomial ideals.  It easily follows
that these Young diagrams are different.  Thus $TRT'\Rightarrow TrT'$.
\end{proof} 

With these propositions, the proof of Proposition \ref{samecomp} is
complete.

\section{Sums of Young diagrams, products of multiple structures}
\label{SoYd,poms}
There is a simple operation on two-dimensional Young diagrams that we are going to need, see
Chapter 7 of Nakajima's book \cite{Nak99}.  Given two partitions
$\lambda=(\lambda_0\geq\lambda_1\geq \dots \geq \lambda_k\geq 0)$ and
$\mu=(\mu_0\geq\mu_1\geq\dots\geq\mu_k\geq 0)$  (we extend the
shortest partition by adding zeroes, if necessary, so that they have
the same length), we can form their partswise sum $(\lambda_1+\mu_1\geq
\dots\geq\lambda_k+\mu_k)$. The result is again a partition.  The
operation on Young diagrams is to add the number of boxes in each
row.  For instance $(4,4,3,2)+(3,3,1)=(4,4,3,2)+(3,3,1,0)=(7,7,4,2)$
or

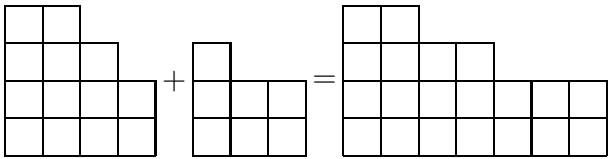
\begin{picture}(4,1)

    \put(0,0){\line(1,0){0.8}}
    \put(0,0.2){\line(1,0){0.8}}
    \put(0,0.4){\line(1,0){0.8}}
    \put(0,0.6){\line(1,0){0.6}}
    \put(0,0.8){\line(1,0){0.4}}

    \put(0,0){\line(0,1){0.8}}
    \put(0.2,0){\line(0,1){0.8}}
    \put(0.4,0){\line(0,1){0.8}}
    \put(0.6,0){\line(0,1){0.6}}
    \put(0.8,0){\line(0,1){0.4}}

    \put(0.9,0.4){\makebox(0,0){$+$}}

    \put(1,0){\line(1,0){0.6}}
    \put(1,0.2){\line(1,0){0.6}}
    \put(1,0.4){\line(1,0){0.6}}
    \put(1,0.6){\line(1,0){0.2}}

    \put(1,0){\line(0,1){0.6}}
    \put(1.2,0){\line(0,1){0.6}}
    \put(1.4,0){\line(0,1){0.4}}
    \put(1.6,0){\line(0,1){0.4}}

    \put(1.7,0.4){\makebox(0,0){$=$}}

    \put(1.8,0){\line(1,0){1.4}}
    \put(1.8,0.2){\line(1,0){1.4}}
    \put(1.8,0.4){\line(1,0){1.4}}
    \put(1.8,0.6){\line(1,0){0.8}}
    \put(1.8,0.8){\line(1,0){0.4}}

    \put(1.8,0){\line(0,1){0.8}}
    \put(2.0,0){\line(0,1){0.8}}
    \put(2.2,0){\line(0,1){0.8}}
    \put(2.4,0){\line(0,1){0.6}}
    \put(2.6,0){\line(0,1){0.6}}
    \put(2.8,0){\line(0,1){0.4}}
    \put(3.0,0){\line(0,1){0.4}}
    \put(3.2,0){\line(0,1){0.4}}
 \end{picture}\\
\noindent
This operation has several good properties:

\begin{prop}
The operation described above is associative, commutative and has unit
given by the empty diagram.  Thus the set of all Young diagrams
becomes an abelian monoid.  This monoid is generated
by Young diagrams with only one column.
\end{prop}
\begin{proof}
The fact that we get an abelian monoid is trivial.  The last statement
follows since any Young diagram is the sum of its columns (or any
partition is the sum of partitions of the form $(1,1,\dots,1)$).
\end{proof}

Based on this operation, we will also construct a three-dimensional
Young diagram as follows: \\

\begin{construction}
\label{konstruksjonav3Ddiag}
Given two partitions $\lambda$ and $\mu$, we make a three-dimensional
Young diagram with boxes
\[\{(i,j,k)|(i,j)\in \mbox{the diagram of }\mu,\mbox{ }k< \lambda_j\}\]
\end{construction}

\begin{bem}
This is {\em not} a commutative construction.  It depends on the order
we write the two diagrams in.
\end{bem}

This can be visualized diagrammatically: consider the skew diagram
given by the right hand side minus the right member of the left side
in the defining sum of two Young diagrams (where $\lambda$ and $\mu$
are added in that order):

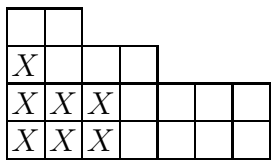
\begin{picture}(4,1)
    \put(0,0){\line(1,0){1.4}}
    \put(0,0.2){\line(1,0){1.4}}
    \put(0,0.4){\line(1,0){1.4}}
    \put(0,0.6){\line(1,0){0.8}}
    \put(0,0.8){\line(1,0){0.4}}

    \put(0,0){\line(0,1){0.8}}
    \put(0.2,0){\line(0,1){0.8}}
    \put(0.4,0){\line(0,1){0.8}}
    \put(0.6,0){\line(0,1){0.6}}
    \put(0.8,0){\line(0,1){0.6}}
    \put(1.0,0){\line(0,1){0.4}}
    \put(1.2,0){\line(0,1){0.4}}
    \put(1.4,0){\line(0,1){0.4}}

    \put(0.1,0.04){\makebox(0,0)[b]{$X$}}
    \put(0.1,0.24){\makebox(0,0)[b]{$X$}}
    \put(0.1,0.44){\makebox(0,0)[b]{$X$}}
    \put(0.3,0.04){\makebox(0,0)[b]{$X$}}
    \put(0.3,0.24){\makebox(0,0)[b]{$X$}}
    \put(0.5,0.04){\makebox(0,0)[b]{$X$}}
    \put(0.5,0.24){\makebox(0,0)[b]{$X$}}

 \end{picture}\\
\noindent
(this corresponds to the sum of $\lambda=(4,4,3,2)$ and $\mu=(3,3,1)$
as above).
The first part consists of those $X$s with at least one empty box to
their right (a {\em part} is a subset of the three dimensional diagram
given by all boxes with a given third coordinate).  This condition
means that $0<\lambda_j$.  The second part
consists of those $X$s that have at least two empty boxes to their
right, meaning $1<\lambda_j$, and so on.  Here the parts are counted
from the first part with third coordinate $0$, to the second part with
third coordinate $1$ and so on.  Thus the
two-dimensional parts of the three-dimensional Young diagram defined by
the above equation are

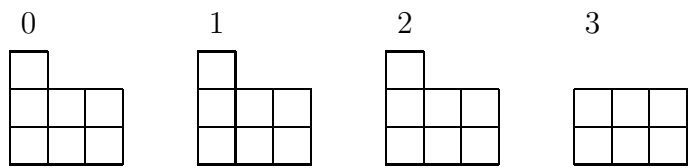
\begin{picture}(4,1)
    \put(0.1,0.7){\makebox(0,0)[b]{$0$}}

    \put(0,0){\line(1,0){0.6}}
    \put(0,0.2){\line(1,0){0.6}}
    \put(0,0.4){\line(1,0){0.6}}
    \put(0,0.6){\line(1,0){0.2}}

    \put(0,0){\line(0,1){0.6}}
    \put(0.2,0){\line(0,1){0.6}}
    \put(0.4,0){\line(0,1){0.4}}
    \put(0.6,0){\line(0,1){0.4}}

   \put(1.1,0.7){\makebox(0,0)[b]{$1$}}

    \put(1,0){\line(1,0){0.6}}
    \put(1,0.2){\line(1,0){0.6}}
    \put(1,0.4){\line(1,0){0.6}}
    \put(1,0.6){\line(1,0){0.2}}

    \put(1,0){\line(0,1){0.6}}
    \put(1.2,0){\line(0,1){0.6}}
    \put(1.4,0){\line(0,1){0.4}}
    \put(1.6,0){\line(0,1){0.4}}

    \put(2.1,0.7){\makebox(0,0)[b]{$2$}}

    \put(2,0){\line(1,0){0.6}}
    \put(2,0.2){\line(1,0){0.6}}
    \put(2,0.4){\line(1,0){0.6}}
    \put(2,0.6){\line(1,0){0.2}}

    \put(2,0){\line(0,1){0.6}}
    \put(2.2,0){\line(0,1){0.6}}
    \put(2.4,0){\line(0,1){0.4}}
    \put(2.6,0){\line(0,1){0.4}}

    \put(3.1,0.7){\makebox(0,0)[b]{$3$}}

    \put(3,0){\line(1,0){0.6}}
    \put(3,0.2){\line(1,0){0.6}}
    \put(3,0.4){\line(1,0){0.6}}

    \put(3,0){\line(0,1){0.4}}
    \put(3.2,0){\line(0,1){0.4}}
    \put(3.4,0){\line(0,1){0.4}}
    \put(3.6,0){\line(0,1){0.4}}

\end{picture}\\
In picture number $i$, the third space coordinate is $i$.
Consider for example the $X$ in position $(1,1)$.  Then $(1,1,k)$ is
in the three-dimensional diagram if and only if $k<\lambda_1=4$, and
$4$ is also the number of empty boxes to the right of this $X$.\\








  



\begin{bem}
If we add the two partitions in the opposite order, and carry through
this construction, the result will be the mirror image of the first
diagram; see Corollary \ref{mirrororder}.
\end{bem}

\begin{bem}
This construction will always give a three-dimensional Young
diagram.  On the other hand, if the skew diagram is not made from
two Young diagrams as above, there is no reason why the result
should be a three-dimensional Young diagram.
\end{bem}

\begin{bem}
If the two partitions we add are
$\lambda=(\lambda_0\geq\lambda_1\geq \dots \geq \lambda_k\geq 0)$ and
$\mu=(\mu_0\geq\mu_1\geq\dots\geq\mu_k\geq 0)$, then the total number
of boxes in the three-dimensional diagram is $\sum_l \lambda_l\mu_l$ since each box
in the diagram corresponding to $\mu$ in the $l$th row is repeated
$\lambda_l$ times.
\end{bem}

Now we consider the operation on monomial ideals
corresponding to this operation on Young diagrams.  For points in the plane, this was
considered by Nakajima, see \cite{Nak99} Chapter 7.  We give an
example from this theory:

\begin{example}
\label{n=0eksempel}
We want to consider two points lying on the line $y=0$ in the
projective plane $\PP^2=\Proj k[x,y,z]$.  Removing the line at
infinity $z=0$, we can coordinatize
this line with the remaining coordinate $x$, and let $x_1,x_2$ be the
points in question.  Let there be given two multiple structures on
these points, with ideals  
$J_1=((x-x_1)^4,(x-x_1)^3y^2,(x-x_1)^2y^3,y^4)$ and
$J_2=((x-x_2)^3,(x-x_2)y^2,y^3)$.  Then the ideal of the union of the
schemes is given by
\[(y^4,y^3(x-x_1)^2,y^2(x-x_1)^3(x-x_2),(x-x_1)^4(x-x_2)^3)\]
\noindent
and the ideal of the special fiber, as $x_2$ tends to $x_1$, becomes
(after the simplification $x_1=0$) 
\[(y^4,y^3x^2,y^2x^4,x^7).\]
The sum of the Young diagrams of the two original ideals is equal to
the Young diagram of the special fibre:

\begin{picture}(4,1)

    \put(0,0){\line(1,0){0.8}}
    \put(0,0.2){\line(1,0){0.8}}
    \put(0,0.4){\line(1,0){0.8}}
    \put(0,0.6){\line(1,0){0.6}}
    \put(0,0.8){\line(1,0){0.4}}

    \put(0,0){\line(0,1){0.8}}
    \put(0.2,0){\line(0,1){0.8}}
    \put(0.4,0){\line(0,1){0.8}}
    \put(0.6,0){\line(0,1){0.6}}
    \put(0.8,0){\line(0,1){0.4}}

    \put(0.9,0.4){\makebox(0,0){$+$}}

    \put(1,0){\line(1,0){0.6}}
    \put(1,0.2){\line(1,0){0.6}}
    \put(1,0.4){\line(1,0){0.6}}
    \put(1,0.6){\line(1,0){0.2}}

    \put(1,0){\line(0,1){0.6}}
    \put(1.2,0){\line(0,1){0.6}}
    \put(1.4,0){\line(0,1){0.4}}
    \put(1.6,0){\line(0,1){0.4}}

    \put(1.7,0.4){\makebox(0,0){$=$}}

    \put(1.8,0){\line(1,0){1.4}}
    \put(1.8,0.2){\line(1,0){1.4}}
    \put(1.8,0.4){\line(1,0){1.4}}
    \put(1.8,0.6){\line(1,0){0.8}}
    \put(1.8,0.8){\line(1,0){0.4}}

    \put(1.8,0){\line(0,1){0.8}}
    \put(2.0,0){\line(0,1){0.8}}
    \put(2.2,0){\line(0,1){0.8}}
    \put(2.4,0){\line(0,1){0.6}}
    \put(2.6,0){\line(0,1){0.6}}
    \put(2.8,0){\line(0,1){0.4}}
    \put(3.0,0){\line(0,1){0.4}}
    \put(3.2,0){\line(0,1){0.4}}

    \put(3.3,0.04){\makebox(0,0)[b]{$x^7$}}
    \put(2.8,0.44){\makebox(0,0)[b]{$x^4y^2$}}
    \put(2.4,0.64){\makebox(0,0)[b]{$x^2y^3$}}
    \put(1.9,0.84){\makebox(0,0)[b]{$y^4$}}
 \end{picture}

The flatness of the deformation can be deduced directly from the
Young diagrams:  the Hilbert polynomial of the multiple point is equal
to the number of boxes in the diagram, and this is ``additive''.\\

Note that the operation is not canonical in the sense that we have
chosen coordinates in the plane; the multiple structures are required
to be ``monomial in the coordinate y and the direction x''.

\end{example}

  Now it is time to describe the changes that must be made if we want
  to extend the results quoted above to higher dimensions.  It turns
  out that there are very few changes, although the arguments become
  harder.  We let $X=\PP^n\subset \PP^{n+2}$.

\begin{prop}
Consider $X=\PP^n\subset\PP^{n+2}$ as contained in a fixed hyperplane
$H=\PP^{n+1}\subset \PP^{n+2}$.  The ideals are $I_X=(x,y)$ and
$I_H=(y)$.  Let $Z=\PP^n$ be
another linear subspace, also contained in $H$, with
ideal $I_Z=(y,z)$.  Then, given Cohen-Macaulay monomial ideals in
$(x,y),(z,y)$ (with these ideals as radicals), let
$z\rightarrow x$.  The family given by the union of the two multiple schemes
on the two linear subspaces $X$ and $Z_t$, where $Z_t$ is determined
by $Z$ under the substitution $z\mapsto tz+(1-t)x$, is a flat
family, and the special fiber ($t=0$) is given by adding the Young
diagrams. 
\end{prop}
\begin{proof}
We can repeat the calculation from Example \ref{n=0eksempel}, and get
the family we need.  The difficulty lies in showing that this family
is flat.  For points, this was simple:  the Hilbert polynomial just
counts the number of boxes in the Young diagram.  In this more
general case the Hilbert polynomial is given by Corollary
\ref{computethehilbpol},
and it is obvious from that description that when we sum Young
diagrams, we don't get the sum of the Hilbert polynomials.  But the
difference is easily calculated.\\

{\em Claim 1.}  Let $T'+T''=T$ be a sum of Young diagrams.  Then the
difference between the Hilbert polynomial of the multiple structure
corresponding
to $T$ and the sum of the Hilbert polynomials of the multiple structures
corresponding to $T'$ and $T''$ is given by the Hilbert polynomial of
the multiple structure corresponding to the three-dimensional
Young diagram from Construction \ref{konstruksjonav3Ddiag}.\\

We also need to give a geometric description of the three-dimensional Young
diagram:\\

{\em Claim 2.} The three-dimensional Young diagram is the Young diagram of the
intersection of the two multiple structures under consideration.\\

Given these two claims, the proposition follows from the general
equality
\[\Hilb(V\cup W)=\Hilb(V)+\Hilb(W)-\Hilb(V\cap W).\]

{\em Proof of Claim 1.} Consider the following diagram, associated to
the addition $T'+T''=T$:

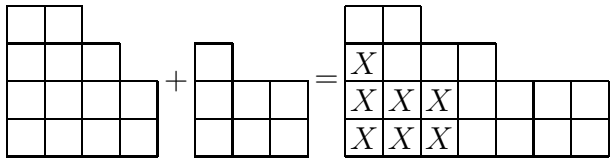
\begin{picture}(4,1)

    \put(0,0){\line(1,0){0.8}}
    \put(0,0.2){\line(1,0){0.8}}
    \put(0,0.4){\line(1,0){0.8}}
    \put(0,0.6){\line(1,0){0.6}}
    \put(0,0.8){\line(1,0){0.4}}

    \put(0,0){\line(0,1){0.8}}
    \put(0.2,0){\line(0,1){0.8}}
    \put(0.4,0){\line(0,1){0.8}}
    \put(0.6,0){\line(0,1){0.6}}
    \put(0.8,0){\line(0,1){0.4}}

    \put(0.9,0.4){\makebox(0,0){$+$}}

    \put(1,0){\line(1,0){0.6}}
    \put(1,0.2){\line(1,0){0.6}}
    \put(1,0.4){\line(1,0){0.6}}
    \put(1,0.6){\line(1,0){0.2}}

    \put(1,0){\line(0,1){0.6}}
    \put(1.2,0){\line(0,1){0.6}}
    \put(1.4,0){\line(0,1){0.4}}
    \put(1.6,0){\line(0,1){0.4}}

    \put(1.7,0.4){\makebox(0,0){$=$}}

    \put(1.8,0){\line(1,0){1.4}}
    \put(1.8,0.2){\line(1,0){1.4}}
    \put(1.8,0.4){\line(1,0){1.4}}
    \put(1.8,0.6){\line(1,0){0.8}}
    \put(1.8,0.8){\line(1,0){0.4}}

    \put(1.8,0){\line(0,1){0.8}}
    \put(2.0,0){\line(0,1){0.8}}
    \put(2.2,0){\line(0,1){0.8}}
    \put(2.4,0){\line(0,1){0.6}}
    \put(2.6,0){\line(0,1){0.6}}
    \put(2.8,0){\line(0,1){0.4}}
    \put(3.0,0){\line(0,1){0.4}}
    \put(3.2,0){\line(0,1){0.4}}

    \put(1.9,0.04){\makebox(0,0)[b]{$X$}}
    \put(1.9,0.24){\makebox(0,0)[b]{$X$}}
    \put(1.9,0.44){\makebox(0,0)[b]{$X$}}
    \put(2.1,0.04){\makebox(0,0)[b]{$X$}}
    \put(2.1,0.24){\makebox(0,0)[b]{$X$}}
    \put(2.3,0.04){\makebox(0,0)[b]{$X$}}
    \put(2.3,0.24){\makebox(0,0)[b]{$X$}}

 \end{picture}\\
\noindent
Here the $X$s in $T$ correspond to the boxes in $T''$.
For each $X$ on the right hand side we get a contribution to the Hilbert
polynomial that is also found on the left hand side in the second
member.  Thus
the difference in Hilbert polynomials must be found using the remaining
boxes on the right hand side and the first member on the left hand side.\\

Consider any other box, and the corresponding box in $T'$.

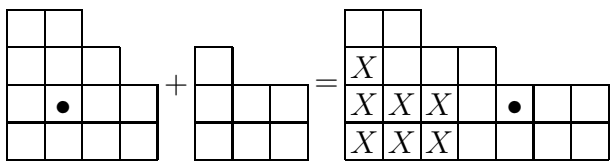
\begin{picture}(4,1)

    \put(0,0){\line(1,0){0.8}}
    \put(0,0.2){\line(1,0){0.8}}
    \put(0,0.4){\line(1,0){0.8}}
    \put(0,0.6){\line(1,0){0.6}}
    \put(0,0.8){\line(1,0){0.4}}

    \put(0,0){\line(0,1){0.8}}
    \put(0.2,0){\line(0,1){0.8}}
    \put(0.4,0){\line(0,1){0.8}}
    \put(0.6,0){\line(0,1){0.6}}
    \put(0.8,0){\line(0,1){0.4}}

    \put(0.3,0.24){\makebox(0,0)[b]{$\bullet$}}

    \put(0.9,0.4){\makebox(0,0){$+$}}

    \put(1,0){\line(1,0){0.6}}
    \put(1,0.2){\line(1,0){0.6}}
    \put(1,0.4){\line(1,0){0.6}}
    \put(1,0.6){\line(1,0){0.2}}

    \put(1,0){\line(0,1){0.6}}
    \put(1.2,0){\line(0,1){0.6}}
    \put(1.4,0){\line(0,1){0.4}}
    \put(1.6,0){\line(0,1){0.4}}

    \put(1.7,0.4){\makebox(0,0){$=$}}

    \put(1.8,0){\line(1,0){1.4}}
    \put(1.8,0.2){\line(1,0){1.4}}
    \put(1.8,0.4){\line(1,0){1.4}}
    \put(1.8,0.6){\line(1,0){0.8}}
    \put(1.8,0.8){\line(1,0){0.4}}

    \put(1.8,0){\line(0,1){0.8}}
    \put(2.0,0){\line(0,1){0.8}}
    \put(2.2,0){\line(0,1){0.8}}
    \put(2.4,0){\line(0,1){0.6}}
    \put(2.6,0){\line(0,1){0.6}}
    \put(2.8,0){\line(0,1){0.4}}
    \put(3.0,0){\line(0,1){0.4}}
    \put(3.2,0){\line(0,1){0.4}}

    \put(1.9,0.04){\makebox(0,0)[b]{$X$}}
    \put(1.9,0.24){\makebox(0,0)[b]{$X$}}
    \put(1.9,0.44){\makebox(0,0)[b]{$X$}}
    \put(2.1,0.04){\makebox(0,0)[b]{$X$}}
    \put(2.1,0.24){\makebox(0,0)[b]{$X$}}
    \put(2.3,0.04){\makebox(0,0)[b]{$X$}}
    \put(2.3,0.24){\makebox(0,0)[b]{$X$}}
    \put(2.7,0.24){\makebox(0,0)[b]{$\bullet$}}
 \end{picture}

The difference in Hilbert polynomials is
$b_i-b_j$ where the box on the left has weight $i$, the box on the
right has weight $j$.  Note that $j-i=$ number of $X$s in that row.
Now there is a fundamental relation of binomial coefficients

\[b_i-b_j=\binom{n+d-i}{n}-\binom{n+d-j}{n}
=\sum_{k=j}^{i-1}\binom{n-1+d-k}{n-1}.\]

For each $X$ in the row, the marked box determines a box in the
three-dimensional diagram above $X$, with third coordinate equal to
the difference between the weight of the marked box and the weight of
the rightmost box with an $X$.  For each $k$ between $j$ and $i-1$ this
gives a box in the three-dimensional diagram with weight $k$.  This
explains the terms of the binomial relation.  For the example marked
above, we have a box in position $(4,1)$ on the right, and $(1,1)$ on
the left.  There are three $X$s in this row, and the three boxes in
the three-dimensional diagram coming from this marked box have
coordinates $(0,1,1),(1,1,1)$ and $(2,1,1)$.\\

The claim follows by letting the marked box run through the skew
diagram (the boxes without $X$s).\\

{\em Proof of claim 2.}  
A box $(i,j,k)$, with corresponding monomial $x^iy^jz^k$, is in
the three-dimensional diagram corresponding to the sum of $\lambda$
and $\mu$ if $(i,j)$ is in the diagram $\mu$
(meaning that the monomial $x^iy^j$ is not in
the ideal corresponding to $\mu$) and if $k<\lambda_j$ (meaning that
$y^jz^k$ is not in the ideal corresponding to $\lambda$.)  This proves
the claim.
\end{proof}

\begin{cor}
The Young diagram of the intersection of two Cohen-Macaulay monomial multiple
structures whose supports are linear subspaces of
codimension two, contained in a common hyperplane, is given by the
construction of a three-dimensional Young diagram from a pair of
two-dimensional Young diagrams.
\end{cor}

\begin{cor}\label{mirrororder}
The operation of constructing a three-dimensional Young diagram from a
pair of two-dimensional Young diagrams does not depend on the order of
the two, up to permutation of the axes.
\end{cor}

\begin{bem}
The permutation of the axes mentioned in the second corollary is given
by permutation of the summands when adding the two Young diagrams.  In
the first corollary, this has the same effect; we will get the two
possible choices of three-dimensional Young diagram by deciding what
order we will write the intersection in.
\end{bem}

\bibliography{ref}
\bibliographystyle{plain}

\end{document}